\documentclass[a4paper,12pt]{article}
\usepackage{amssymb,amsmath}
\newcommand{\Q}{\mathbb{Q}}

\newcommand{\N}{\mathbb{N}}
\newcommand{\Proj}{\mathbb{P}}
\renewcommand{\b}[1]{{\bf #1}}
\newcommand{\cl}[1]{{\cal #1}}
\newcommand{\Z}{\mathbb{Z}}

\newcommand{\ep}{\varepsilon}
\newcommand{\y}{\b{y}}
\newcommand{\x}{\b{x}}

\newtheorem{theorem}{Theorem}
\newtheorem{lemma}{Lemma}

\newcommand{\la}{\lambda_1,\lambda_2}
\newcommand{\ta}{\tau_1,\tau_2}

\begin{document}
\title{Square-Free Values of $n^2+1$}
\author{D.R. Heath-Brown\\Mathematical Institute, Oxford}
\date{}
\maketitle

\begin{flushright} 
To Professor Andrzej Schinzel\\
In Celebration of his Seventy-Fifth Birthday
\end{flushright}

\section{Introduction}
Let $\cl{N}(x)$ denote the number of positive integers $n\le x$ for
which $n^2+1$ is square-free.  It was shown in 1931 by Estermann \cite{Est} that
\[\cl{N}(x)=c_0x+O(x^{2/3}\log x)\]
for $x\ge 2$, where
\[c_0=\frac{1}{2}\prod_{p\equiv 1\bmod{4}}(1-2p^{-2}).\]
Estermann's argument is very simple, but despite the passage of 80 years the
exponent $2/3$ appearing above has never been improved.  The aim of
the present paper is to establish the following result.
\begin{theorem}
We have
\[\cl{N}(x)=c_0x+O_{\ep}(x^{7/12+\ep})\]
for any fixed $\ep>0$.
\end{theorem}

It is easy to construct intervals $(x,x+c\log x]$ with a small
positive constant $c$, such that $n^2+1$ has a non-trivial square
factor for every $n$ in the interval.  This shows that the error term
in our theorem is $\Omega(\log x)$.  However we know of no better
result of this type, and it is unclear what one should conjecture.
With the much simpler problem of the number of square-free integers
$n\le x$ one has an easy error term $O(x^{1/2})$, but any reduction in
the exponent $1/2$ would appear to require a quasi Riemann Hypothesis.
Thus it seems unlikely that we could reduce the exponent in our
theorem below $1/2$ without a radically new idea.

The key point in our treatment will be to give good upper bounds for
the frequency of solutions to the Diophantine equation $d^2e=n^2+1$.
Analysing this over $\Q(i)$ we are led to study the condition
$2x_1x_2y_1+(x_1^2-x_2^2)y_2=1$, which we may interpret as saying that
the point $(s,t)=(x_1/x_2,y_1/y_2)$ lies close to the curve
$t=(1-s^2)/(2s)$. In order to study this we we will use a variant of 
the ``Determinant Method'',
developed from the author's papers \cite{Ann}, \cite{HB}.  \medskip

The author was introduced to this problem by Dr Tim Browning.  His
contributions to the resulting discussions, and his careful
proof-reading of the present paper, are gratefully acknowledged.

\section{Preliminaries}\label{prelim}

For the proof it will clearly suffice to show that
\[\cl{N}(2x)-\cl{N}(x)=c_0x+O_{\ep}(x^{7/12+\ep})\]
The argument begins by observing that for $x\ge 1$ and $1\le D\le x$ we have
\begin{eqnarray}\label{1}
\cl{N}(2x)-\cl{N}(x)&=&\sum_{x<n\le 2x}\mu^2(n^2+1)\nonumber\\
&=&\sum_{x<n\le 2x}\;\sum_{d^2\mid n^2+1}\mu(d)\nonumber\\
&=&\sum_{d\le 4x}\mu(d)\#\{x<n\le 2x: d^2\mid n^2+1\}\nonumber\\
&=&\sum_{d\le D}\mu(d)\#\{x<n\le 2x: d^2\mid n^2+1\}\nonumber\\
&&\hspace{1cm}\mbox{}+O(\sum_{D<d\le 4x}\#\{x<n\le 2x: d^2\mid n^2+1\}).
\end{eqnarray}
For $d\le D$ we write
\[\rho(d)=\rho=\#\{m\bmod{d^2}:d^2\mid m^2+1\},\]
and we take $m_1,\ldots,m_{\rho}$ to be a corresponding set of
admissible values for $m$.  Then
\begin{eqnarray}\label{dr1}
\#\{x<n\le 2x: d^2\mid n^2+1\}&=&
\sum_{j=1}^{\rho}\#\{x<n\le 2x: n\equiv m_j\bmod{d^2}\}\nonumber\\
&=&\sum_{j=1}^{\rho}\left(\frac{x}{d^2}+O(1)\right)\nonumber\\
&=&x\frac{\rho(d)}{d^2}+O(\rho(d)).
\end{eqnarray}
Thus terms with $d\le D$ contribute to (1) a total
\[x\sum_{d\le D}\mu(d)\rho(d)d^{-2}+O(\sum_{d\le D}\rho(d)).\]
The function $\rho(d)$ is multiplicative, with $\rho(p)=2$ for
$p\equiv 1\bmod{4}$ and $\rho(p)=0$ otherwise.  Thus $\rho(d)$ is
bounded by the familiar $r(d)$ function which counts representations
as sums of two squares.  We therefore see that
\[\sum_{E<d\le 2E}\rho(d)\ll E\]
for any integer $E$, whence
\[\sum_{d>D}|\mu(d)\rho(d)d^{-2}|\ll D^{-1}\;\;\;\mbox{and}\;\;\;
\sum_{d\le D}\rho(d)\ll D.\]
The contribution to (1) corresponding to values $d\le D$ is therefore
\[x\sum_{d=1}^{\infty}\frac{\mu(d)\rho(d)}{d^2}+O(xD^{-1})+O(D).\]
Since
\[\sum_{d=1}^{\infty}\frac{\mu(d)\rho(d)}{d^2}=\prod_p(1-\rho(p)p^{-2})\]
we see that this produces the main term in our theorem.  We will
minimize the other error terms by choosing $D=x^{1/2}$.

To handle the larger values of $d$ we consider dyadic ranges $E/2<d\le
E$, and write
\[\cl{M}(E,F)=\#\{(e,f,n)\in\N^3: E/2<e\le E,\, F/2<f\le F,\,
e^2f=n^2+1\}.\]
Then the range $d>D$ contributes to (1) a total
\[\ll\sum_{E\gg D}\;\max_{x^2E^{-2}\ll F\ll x^2E^{-2}}\cl{M}(E,F),\]
where the summation for $E$ runs over powers of 2.
Thus our problem reduces to one of estimating $\cl{M}(E,F)$
efficiently.  Heuristically one might expect that $e^2f-1$ is a square
``with probability'' of order $(e^2f)^{-1/2}$.  This leads one to
conjecture that the true order of magnitude for $\cl{M}(E,F)$ might be
about $F^{1/2}$. For his proof, Estermann showed that
\begin{equation}\label{E}
\#\{(e,n)\in\N^2: E/2<e\le E,\, e^2f=n^2+1\}\ll \log x,
\end{equation}
whence $\cl{M}(E,F)\ll F\log x$.
One easily sees how this leads to the error term $O(x^{2/3}\log
x)$.  We will need sharper bounds, but we note that Estermann's
estimate shows that the range $F\le x^{1/2}$ yields a satisfactory
contribution.  Since we have taken $D=x^{1/2}$ we may therefore assume 
in what follows that
$x^{1/2}\ll E\ll x^{3/4}$ and $x^{1/2}\ll F\ll x$.

\section{The Determinant Method}
We begin our analysis of $\cl{M}(E,F)$
by using the unique factorization property for
$\Z[i]$. This shows that if $e^2f=n^2+1$ then there are integers
$x_1,x_2,y_1,y_2$ for which 
\[e=x_1^2+x_2^2,\;\;\;f=y_1^2+y_2^2\;\;\;\mbox{and}
\;\;\;(x_1+ix_2)^2(y_1+iy_2)=n+i.\]
It follows on taking the imaginary part
that $2x_1x_2y_1+(x_1^2-x_2^2)y_2=1$.  If $|x_1|>|x_2|$ we will swap
$x_1$ and $x_2$, and change the sign of $y_2$.  Hence we may suppose,
without loss of generality, that $|x_1|\leq |x_2|$, and hence that
$|x_1|\le E^{1/2}$ and $E^{1/2}\ll |x_2|\le E^{1/2}$.
We observe that
\[\max\{|2x_1x_2|\,,\,|x_1^2-x_2^2|\}\gg x_1^2+x_2^2\gg E.\]
Thus if we write $q_1(x_1,x_2)=2x_1x_2$ or $x_1^2-x_2^2$ as
appropriate, and take $q_2$ to be the alternative quadratic form, we
may assume that $E\ll |q_1(x_1,x_2)|\ll E$ and $q_2(x_1,x_2)\ll E$.
Then, labelling $y_1,y_2$ either as $z_1,z_2$ or as $z_2,z_1$ we will
have 
\begin{equation}\label{2}
q_1(x_1,x_2)z_1+q_2(x_1,x_2)z_2=1, 
\end{equation}
whence
\[E|z_1|\ll |q_1(x_1,x_2)z_1|\ll 1+|q_2(x_1,x_2)z_2|\ll 1+E|z_2|.\]
Since we cannot have $z_1=z_2=0$ we deduce that $|z_1|\ll |z_2|$.  
Then, since $F\ll z_1^2+z_2^2\le F$ we see that 
$|z_1|\le F^{1/2}$ and $F^{1/2}\ll |z_2|\le F^{1/2}$.

We now deduce from (\ref{2}) that if $s=x_1/x_2$ and $t=z_1/z_2$ then
\[t=-\frac{q_2(s,1)}{q_1(s,1)}+O(E^{-1}F^{-1/2})
=-\frac{q_2(s,1)}{q_1(s,1)}+O(x^{-1}).\]
Thus if we write $\phi(s)=-q_2(s,1)/q_1(s,1)$ then the point $(s,t)$
lies close to the curve $\phi(S)=T$. Our task is therefore to estimate
the number of rational points $(s,t)$ with $s,t\ll 1$ lying within
$O(x^{-1})$ of the curve $\phi(S)=T$, and for which the ``heights'' of
$s$ and $t$ are at most $E^{1/2}$ and $F^{1/2}$ respectively.

The situation here is similar to that in the author's paper \cite{HB}.
We shall use a real-variable version of the ``determinant method'',
but there is an important difference, in that the variety given by the
equation 
\[q_1(x_1,x_2)z_1+q_2(x_1,x_2)z_2=0\]
lies naturally in
$\Proj^1\times \Proj^1$, rather than in $\Proj^2$.  Indeed this makes
our situation correspond exactly to that considered by Huxley
\cite{Hux1}, \cite{Hux2}.  Unfortunately Huxley's bounds, which were
obtained for general plane curves, are not strong enough for our
application. In particular, he focuses on the case in which $\phi$ is
not a rational function.

Following the method from the author's work \cite[\S 2]{HB} we choose
an integer parameter $M\in [x^{1/2},x]$ and split
the available range for $s$ into $O(M)$ subintervals $I=(s_0,s_0+M^{-1}]$.
We then investigate the number of solutions in which $s$ belongs to a
particular interval $I$.  If we 
write $s=s_0+u$ we find from Taylor's Theorem that
$\phi(s)=\phi(s_0)+u\phi'(s_0)+O(x^{-1})$. Hence if we set 
$v=t-\phi(s_0)-u\phi'(s_0)$ we will have $s=s_0+u, 
t=\phi(s_0)+u\phi'(s_0)+v$ with $v\ll
x^{-1}$.  We now label all the solutions corresponding to the interval
$I$ as $(s_1,t_1),\ldots,(s_J,t_J)$, say.
We proceed to choose positive integers $K,L$ and to
label the monomials $s^kt^l$ for $k\le K$, $l\le L$ as
$m_1(s,t),\ldots,m_H(s,t)$, where $H=(K+1)(L+1)$.  The determinant
method uses the $J\times H$ matrix $\cl{M}$, whose $jh$ entry is
$m_h(s_j,t_j)$.  The aim is to show that the rank of $\cl{M}$ is
strictly less than $H$.  If this can be achieved, one may deduce that there
is a non-zero vector $\mathbf{c}$ with 
\begin{equation}\label{3}
\cl{M}\mathbf{c}=\mathbf{0}.
\end{equation}
This vector $\mathbf{c}$ may be constructed out of appropriate
subdeterminants of $\cl{M}$.  Thus its entries will be rational
numbers with numerators and denominators of size $\ll_{K,L} x^{H(K+L)}$, since
$s_j$ and $t_j$ have numerators and denominators of size $\ll
x^{1/2}\ll x$.
We now observe that the matrix equation (\ref{3}) means that 
there is a polynomial
$C(s,t)$, with coefficients given by the vector $\mathbf{c}$, such
that $C(s_j,t_j)=0$ for all pairs $(s_j,t_j)$.  Multiplying out the
common denominator of the coefficients we may assume that $C$ has
integer coefficients, of size $\ll_{K,L}x^{H^2(K+L)}$.

This is one of the key stages in the proof.  We deduce that all points
$(s,t)$ for which $t$ is close to $\phi(s)$, and for which $s$ lies in an
appropriate short range $I$, actually lie on the curve $C(s,t)=0$.

We now show that $\cl{M}$ does indeed have rank less than $H$, if the
parameter $M$ is suitable chosen.  For this we select any $H\times H$
subdeterminant, $\Delta$ say, from $\cl{M}$, and show that
$\Delta=0$. Without loss of generality we may suppose that $\Delta$
comes from the first $H$ rows of $\cl{M}$.  Since the $j$-th row
contains rationals with a common denominator of $x_{2,j}^Kz_{2,j}^L$ it
is clear that
\[\left(\prod_{j\le H}x_{2,j}^Kz_{2,j}^L\right)\Delta\in\Z.\]
Thus to show that $\Delta=0$ it will suffice to prove that
\begin{equation}\label{add}
\Delta\ll_{K,L}E^{-KH/2}F^{-LH/2}
\end{equation}
with a suitably small implied constant.

When we substitute $s=s_0+u$ and $t=\phi(s_0)+u\phi'(s_0)+v$ the 
monomials $m_j(s,t)$ 
produce polynomials in $u,v$.  Thus $\Delta$ is a generalized van der
Monde determinant.  If $u_j$ and $v_j$ correspond to $s_j$ and $t_j$
then we have $|u_j|\le M^{-1}$ and $|v_j|\le V^{-1}$ for some $V$ of
exact order $x$. An estimate for the size of $\Delta$ is now provided
by Lemma 3 from the author's work \cite{HB}. If we order all
possible monomials $M^{-k}V^{-l}$ in decreasing size as $1=M_0,M_1,\ldots$
then the lemma shows that 
\[\Delta\ll_H\prod_{h=1}^H M_h.\]
If $M_H=W^{-1}$ then $M^{-k}V^{-l}\ge M_H$ if and only if 
\begin{equation}\label{4}
k\log M+l\log V\le\log W.  
\end{equation}
The number of such pairs $k,l$ is 
\[\frac{(\log W)^2}{2(\log M)(\log V)}
+O\left(\frac{\log W}{\log x}\right)+O(1),\]
and since this must equal $H$ we deduce that
\begin{equation}\label{5}
\log W=H^{1/2}\sqrt{2(\log M)(\log V)}+O(\log x).
\end{equation}
Moreover 
\begin{eqnarray*}
\log\prod_{h=1}^H M_h&=&-\sum_{k,l}(k\log M+l\log V)\\
&=&-\frac{(\log W)^3}{3(\log M)(\log V)}+
O\left(\frac{(\log W)^2}{\log x}\right),
\end{eqnarray*}
the sum over $k,l$ being subject to (\ref{4}).  It follows from (\ref{5}) that
\[\log\prod_{h=1}^H M_h=
-H^{3/2}\frac{2\sqrt{2}}{3}\sqrt{(\log M)(\log V)}+O(H\log x),\]
and hence that
\[\log|\Delta|\le O_H(1)
-H^{3/2}\frac{2\sqrt{2}}{3}\sqrt{(\log M)(\log V)}+O(H\log x).\]
This will be sufficient for (\ref{add}) providing that
\[\frac{K}{2}\log E+\frac{L}{2}\log F\le 
(KL)^{1/2}\frac{2\sqrt{2}}{3}\sqrt{(\log M)(\log V)}+
O_{K,L}(1)+O(\log x).\]
In order to use this optimally we will take $K=[L(\log F)/(\log E)]$.
Since our size constraints on $E$ and $F$ imply that $L\ll K\ll L$ it 
then suffices that
\[L\log F\le 
L\frac{2\sqrt{2}}{3}\sqrt{(\log M)(\log V)}\frac{\sqrt{\log F}}{\sqrt{\log E}}
+O_L(1)+O(\log x).\]
Hence if $\delta>0$ is a small positive constant, and 
\[\frac{2\sqrt{2}}{3}\sqrt{(\log M)(\log V)}\frac{\sqrt{\log F}}{\sqrt{\log E}}
\ge (1+\delta)\log F\]
it will be enough to have $L=L(\delta)$ sufficiently large, 
and $x\gg_{\delta}1$.  The condition may be rewritten in the form
\[\log M\ge \tfrac98(1+\delta)^2\frac{(\log E)(\log F)}{\log V}\]
and since $V\gg x$ we may summarize our conclusions as follows.
\begin{lemma}\label{y}
Let $\eta>0$ be given, and suppose $M\in[x^{1/2},x]$ satisfies
\[\log M\ge \tfrac98(1+\eta)\frac{(\log E)(\log F)}{\log x}.\]
Then for any interval $I=[s_0,s_0+M^{-1}]$ there is a corresponding
non-zero integer polynomial $C_I(s,t)$ satisfying
\begin{equation}\label{aux}
C_I(x_1/x_2,z_1/z_2)=0
\end{equation}
for any solution of (\ref{2}) with $x_1/x_2\in I$.
Moreover $C_I$ has total degree $O_{\eta}(1)$, and coefficients of 
size $O_{\eta}(x^{\kappa})$ for some constant $\kappa=\kappa(\eta)$.
\end{lemma}

\section{Counting Solutions of Equations}

While the previous section involved the application of a general
method, the next stage in the proof requires an {\em ad hoc}
argument, to count points which simultaneously satisfy both (\ref{2}) and
(\ref{aux}).  We begin by
showing that it suffices to assume that $C_I$ is absolutely
irreducible. Let $(s,t)$ be a 
rational point satisfying $F(s,t)=0$ for
some monic factor $F$ of $C_I$ which is not defined over $\Q$.  Then
$F^{\sigma}(s,t)=0$ for every conjugate $F^{\sigma}$.  The
number of possible points $s,t$ is then $O_{\eta}(1)$ by B\'ezout's
Theorem. Since $x_1,x_2$ are coprime, and similarly for $z_1,z_2$ we
obtain $O_{\eta}(1)$ solutions this way.  Thus we need only consider
absolutely irreducible factors $F$ of $C_I$ which are defined over
$\Z$. The height of any such factor is again bounded by a power of $x$,
by Gelfond's Lemma (see Bombieri and Gubler \cite[Lemma 1.6.11]{BG}
for example).  Moreover the number of different
factors to consider is $O_{\eta}(1)$.  Thus it suffices to consider
the case in which $F(x_1/x_2,z_1/z_2)=0$ for some absolutely
irreducible polynomial $F$ satisfying the same conditions as $C_I$.

Our next move is to clear the denominators $x_2$
and $z_2$ so as to replace the equation $F(s,t)=0$ by a bi-homogeneous
one \begin{equation}\label{bh}
F(x_1,x_2;z_1,z_2)=0, 
\end{equation}
say. For a given interval $I$ we will have 
\[s_0<s=\frac{x_1}{x_2}\le s_0+M^{-1}.\]
It therefore follows that $|x_1-s_0x_2|\le E^{1/2}M^{-1}$, since
$|x_2|\le E^{1/2}$.  If we let
\[\Lambda=\{\left(E^{-1/2}M(x_1-s_0x_2)\,,\,E^{-1/2}x_2\right):
(x_1,x_2)\in\Z^2\}\] 
then $\Lambda$ is a lattice of determinant $E^{-1}M$, and we are
interested in points $(\alpha_1,\alpha_2)\in\Lambda$ falling in the square 
\[S=\{(\alpha_1,\alpha_2): \max(|\alpha_1|,|\alpha_2|)\le 1\}. \]
Let $\b{g}^{(1)}$ be the shortest non-zero vector in the lattice and
$\b{g}^{(2)}$ the shortest vector not parallel to $\b{g}_1$.  These vectors
will form a basis for $\Lambda$.  Moreover we have 
$\lambda_1\b{g}^{(1)}+\lambda_2\b{g}^{(2)}\in S$ only when $|\lambda_1|\ll
|\b{g}^{(1)}|^{-1}$ and $|\lambda_2|\ll|\b{g}^{(2)}|^{-1}$. These  
constraints may be written in the form
$|\lambda_i|\le L_i$, for appropriate bounds $L_1,L_2$.
Since $|\b{g}^{(2)}|\ge|\b{g}^{(1)}|$ and 
$|\b{g}^{(1)}|.|\b{g}^{(2)}|\ll {\rm det}(\Lambda)=E^{-1}M$ we will have
$L_1\gg L_2$ and $L_1L_2\gg EM^{-1}$.  We now write
\[\b{h}^{(i)}=E^{1/2}(M^{-1}g^{(i)}_1+s_0g^{(i)}_2\,,\,g^{(i)}_2)\]
for  $i=1,2$.  These vectors will then be a basis for $\Z^2$, and if
$\b{x}=\lambda_1\b{h}^{(1)}+\lambda_2\b{h}^{(2)}$ is in the region
given by $|x_1-s_0x_2|\le E^{1/2}M^{-1}$ and $|x_2|\le E^{1/2}$ then
we will have $|\lambda_i|\le L_i$ for $i=1,2$.  This allows us to make
a change of basis, replacing $(x_1,x_2)$ by $(\lambda_1,\lambda_2)$
so that our constraints on $x_1,x_2$ are replaced by the conditions
$|\lambda_i|\le L_i$.

We may argue in exactly the same way for $z_1,z_2$ using the fact 
that 
\[t=z_1/z_2=\phi(s_0)+u\phi'(s_0)u+O(x^{-1})=\phi(s_0)+O(M^{-1}).\]
This allows us to replace the
variables $z_1,z_2$ by $\tau_1,\tau_2$ subject to $|\tau_i|\le T_i$.
Here $T_1\gg T_2$ and $T_1T_2\ll FM^{-1}$.  These substitutions 
convert (\ref{2}) into a new
equation
\begin{equation}\label{G0}
G_0(\la;\ta)=1
\end{equation}
say, where $G_0$ is bi-homogeneous of degree $(2,1)$.  Similarly they
will turn (\ref{bh}) into an equation of the shape
\begin{equation}\label{G1}
G_1(\la;\ta)=0
\end{equation}
where $G_1$ is bi-homogeneous of degree $(a,b)$, say.  Of course it is apparent
from (\ref{G0}) that the vectors $(\la)$ and $(\ta)$ will be primitive.

When $\min(a,b)\ge 2$ we can get a satisfactory bound from the
following general result, which will be proved later, in \S \ref{PC}.
\begin{lemma}\label{cover}
Let $G(x_1,x_2;y_1,y_2)\in\Z[x_1,x_2,y_1,y_2]$ be an absolutely
irreducible bi-homogeneous
polynomial of degree $(a,b)$ with $a,b\ge 1$.  Let $\ep>0$ be
given. Then for any $X\ge 1$ there are 
$O_{a,b,\ep}(X^{2/b+\ep}||G||^{\ep})$ points
$(x_1,x_2,y_1,y_2)\in\Z^4$ satisfying the conditions
\[{\rm g.c.d.}(x_1,x_2)=1,\,\;\;{\rm g.c.d.}(y_1,y_2)=1,\]
\[G(x_1,x_2;y_1,y_2)=0,\;\;\;\mbox{{\rm and}}\;\;\;\max_i|x_i|\le X.\]
\end{lemma}
Notice here that there is no size constraint on $y_1$ or $y_2$. 

When $a\ge 2$ the lemma shows
that (\ref{G1}) has $O_{\eta}(T_1^{1+\eta}x^{\eta})$ solutions
$\lambda_1,\lambda_2$, $\tau_1,\tau_2$.  Each of these corresponds to at
most one solution of (4), and therefore contributes $O(1)$ to $\cl{M}(E,F)$.
Similarly, if $b\ge 2$ that there are $O_{\eta}(L_1^{1+\eta}x^{\eta})$
solutions.  

We next
dispose of the case in which $a$ or $b$ is zero.  For example, if
$a=0$ then (\ref{G1}) specifies a finite number $O_{\eta}(1)$ of 
pairs $\ta$, and for each of these there is a corresponding pair
$(z_1,z_2)$, producing a value of $f=z_1^2+z_2^2$.  Each such $f$
contributes $O(\log x)$ to $\cl{M}(E,F)$ by Estermann's bound
(\ref{E}).  Thus the interval $I$ contributes $O_{\eta}(\log x)$ when
$a=0$.  Similarly if $b=0$ there are $O_{\eta}(1)$ corresponding
values for $e$. As in (\ref{dr1}) each such value $e$ contributes $\ll
\rho(e)(xe^{-2}+1)$ to $\cl{M}(E,F)$.  This is also satisfactory,
since 
\[e\ge E\gg D=x^{1/2}\]
and $\rho(e)\ll_{\eta}x^{\eta}$.  Thus we
have $O_{\eta}(x^{\eta})$ solutions corresponding to $I$ when $\min(a,b)=0$.

When $a=1$ the equation (\ref{G1}) can be written
\[\lambda_1G_{11}(\ta)+\lambda_2G_{12}(\ta)=0.\]
Thus
$\lambda_1=q^{-1}G_{12}(\ta)$, $\lambda_2=-q^{-1}G_{11}(\ta)$, where
$q$ divides $G_{11}(\ta)$ and $G_{12}(\ta)$.  Since $\tau_1$ and
$\tau_2$ are coprime it follows that $q$ divides the resolvent $R$ of
$G_{11}$ and $G_{12}$. This resolvent is non-zero since $G_1$ is
irreducible.  Moreover it is bounded by a power of $x$, whence there
are $O_{\eta}(x^{\eta})$ possible choices for $q$.  (The reader should
recall at this point that the forms $G_{11}$ and $G_{12}$ are determined, up
to $O_{\eta}(1)$ possibilities, by the interval $I$.)  For each
available choice of $q$
we substitute our values for $\la$ into (\ref{G0}) to obtain a Thue
equation $G_3(\ta)=q^2$.  Unfortunately we cannot use the full force of
known results on such equations, since it is possible for $G_3$ to be
a power of a linear form.  None the less there can be at most $O(T_1)$
possible pairs $\ta$.  It follows that we have
at most $O_{\eta}(x^{\eta}T_1)$ solutions in total.  The case $b=1$ is
entirely analogous, leading to a bound $O_{\eta}(x^{\eta}L_1)$.

In summary we have a bound $O_{\eta}(T_1^{1+\eta}x^{\eta})$ on the
number of solutions, in each of the cases $a\ge 2$, or $a=0$, or
$a=1$.  Similarly we have an estimate $O_{\eta}(L_1^{1+\eta}x^{\eta})$
whatever the value of $b$.  We therefore conclude as follows.
\begin{lemma}\label{x}
For any $\eta>0$
the contribution to $\cl{M}(E,F)$ corresponding to a single interval
$I$ is $O_{\eta}(x^{\eta}\min(L_1^{1+\eta},T_1^{1+\eta}))$.
\end{lemma}

\section{Completion of the Proof}

Having fixed $M$ as in Lemma 1 we must now sum up $\min(L_1,T_1)$ for
the various intervals $I$.  In the notation of the previous section,
if $(x_1,x_2)$ corresponds to $\b{g}^{(1)}$ then $L_1(x_1-s_0x_2)\ll
M^{-1}\sqrt{E}$ and $L_1x_2\ll\sqrt{E}$.  If $L_1\gg\sqrt{E}$ we see
that $x_2=0$, and then $x_1=0$, which is impossible. The intervals
$I=(s_0,s_0+M^{-1}]$ will be produced by taking $s_0=x_3M^{-1}$ for
integers $x_3\ll M$.  Thus the number of intervals for which $L\le
L_1\le 2L$ is at most the number of triples $(x_1,x_2,x_3)\in\Z^3$
with ${\rm g.c.d.}(x_1,x_2)=1$, for which
\[x_2x_3=Mx_1+O(L^{-1}\sqrt{E}),\;\;\;
x_2\ll L^{-1}\sqrt{E},\;\;\;\mbox{and}\;\;\;x_3\ll M.\]
We now recall that $L_1\gg L_2$ and that $L_1L_2\gg EM^{-1}$.  Thus
$L\gg E^{1/2}M^{-1/2}$.  Moreover, as noted above, we have $L\ll
E^{1/2}$. In particular, if $M$ is large enough we
can have $x_2x_3=0$ only when $x_1=0$.  Since $x_1$ and $x_2$ are
coprime this case can arise only when $x_2=\pm 1$ and $x_3=0$.
When $x_2x_3\not=0$ the conditions on $x_2$ and $x_3$ imply that 
$x_1\ll L^{-1}\sqrt{E}$, and a divisor function estimate then shows 
that there are
$O_{\eta}(x^{\eta}L^{-1}\sqrt{E})$ pairs $x_2,x_3$ for each value of $x_1$.  
We conclude that there are $O_{\eta}(x^{\eta}L^{-2}E)$ intervals $I$
for which $L_1$ is of order $L$.  Since each interval makes a contribution
$\ll_{\eta} x^{\eta}L_1^{1+\eta}$, by Lemma~\ref{x},  we get a total 
$\ll_{\eta}x^{\eta}L^{-1+\eta}E\ll x^{2\eta}L^{-1}E$, since $L\ll
E^{1/2}\ll x$.  By dyadic subdivision for $L\gg
E^{1/2}M^{-1/2}$ we find that
$\cl{M}(E,F)\ll_{\eta}x^{2\eta}(EM)^{1/2}$.  

We can prove a precisely analogous estimate $\cl{M}(E,F)\ll_{\eta}
x^{2\eta}(FM)^{1/2}$ by considering the number of intervals
$J=(\phi(s_0),\phi(s_0)+O(M^{-1})]$ which produce a value $T_1$ in a
given dyadic range $(T,2T]$.  Here we use the fact that $J\subseteq
(t_3M^{-1},t_3M^{-1}+O(M^{-1})]$ for some integer $t_3$. We also need
to remark that each value of $t_3$ occurs $O(1)$ times, since
$|\phi'(s)|\gg 1$ for the values of $s$ under consideration.  With
these observations the argument then goes through just as before.
We may therefore conclude that
\[\cl{M}(E,F)\ll_{\eta}x^{2\eta}(\min(E,F)M)^{1/2}.\]

It remains to use this result with the value for $M$ coming from
Lemma \ref{y}.  It is convenient to write $E=x^{\psi}$, so that 
$F=x^{2-2\psi+O(1/\log x)}$.  In view of our remarks at the end of \S
2 we have (essentially) $1/2\le\psi\le 3/4$.
We may then employ a value $M$ with
\[\frac{\log M}{\log x}= 
(1+\eta)\max\left\{\frac{9\psi(1-\psi)}{4}\,,\,\frac{1}{2}\right\}
+O((\log x)^{-1}).\]
This value will automatically satisfy $M\in[x^{1/2},x]$ if $\eta$ is
small enough.  It follows that
\begin{eqnarray*}
\frac{\log \cl{M}(E,F)}{\log x}&\le &2\eta+\tfrac12\min(\psi,2-2\psi)+
(1+\eta)\max\left\{\frac{9\psi(1-\psi)}{8}\,,\,\frac{1}{4}\right\}\\
&&\hspace{1cm}\mbox{}+O_{\eta}((\log x)^{-1}).
\end{eqnarray*}
However, since
\[\frac{1}{2}\min(\psi,2-2\psi)+  %
\max\left\{\frac{9\psi(1-\psi)}{8}\,,\,\frac{1}{4}\right\}
\le \frac{7}{12}\] 
for the relevant range of $\psi$,
we deduce that $\cl{M}(E,F)\ll_{\eta} x^{3\eta+7/12}$, and our theorem
then follows.

\section{Lemma \ref{cover}}\label{PC}

Lemma \ref{cover} is closely related to two results of Broberg.  In
\cite[Theorem 1]{br1} Broberg establishes a general result about finite
covers of $\Proj^1$ which, when translated into our notation, would
provide an estimate $O_{G,a,b,\ep} (X^{2/b+\ep})$ of the desired
order, but without any explicit dependence on $G$.  This explicit dependence
can be deduced from a second result of Broberg
\cite{br2}, but this has not been formally published.  We therefore
give a brief sketch of a direct argument independent of these two papers.
This uses the determinant method, for more details of which the reader should 
consult \cite[\S 3]{Ann}.

We will need a crude bound on the size of $\y$.  Let
\[G(\x;\y)=G_0(\x)y_1^b+\ldots+G_b(\x)y_2^b.\]
The form $G_0$ cannot vanish identically since $G$ is irreducible.  Thus
there are $O_a(1)$ primitive integer vectors $\x$ for which $G_0(\x)=0$.  It
is not possible for all the forms $G_i(\x)$ to vanish simultaneously
for a vector $\x\not=\b{0}$, since $G$ is irreducible.  Thus, with
$O_{a,b}(1)$ exceptions, any solution of $G(\x;\y)=0$ has
$y_2|G_0(\x)$, with $G_0(\x)\not=0$.
We may therefore assume that $|y_2|\ll X^a||G||$, and similarly for $y_1$.
It will be convenient to write these bounds in the form
$|y_1|,|y_2|\le Y$, with $Y\ll X^a||G||$.

Our overall plan now is to apply the $p$-adic determinant method.  By
making an invertible integral linear substitution on $\x$ we may
assume that $G_0(1,0)\not=0$. Indeed we can choose the coefficients of
the substitution to be bounded in terms of $a$ alone, so that we may
still assume that $|\x|\ll_a X$.
Suppose we have a parameter $P\gg_{a,b} \log X||G||$.  Then for any
solution $\x,\y$ of $G(\x;\y)=0$ with $|\x|\ll_a X$, $|\y|\ll Y$,
we either have
\[x_2y_2\frac{\partial G(\x;\y)}{\partial y_1}=0\]
or there is a prime $p\in(P,2P]$ not dividing $x_2y_2(\partial
G/\partial y_1)$.  The first case immediately give us an
auxiliary bi-homogeneous form $H(\x;\y)$ not divisible by $G$, at
which our solution also vanishes.  In the alternative case the point
$(x_1/x_2,y_1/y_2)$ lies above a smooth $\mathbb{F}_p$-point on the
curve $G(s,1;t,1)=0$.  We can then expand $y_1/y_2$ as a $p$-adic power
series in $x_1/x_2$, as in Lemma~5 of the author's paper \cite{Ann}.
We then consider the matrix of bi-homogeneous monomials in $\x$ and
$\y$ of degree $(H,b-1)$.  There are $k:=(H+1)b$ such monomials.  The
corresponding $k\times k$ determinant then has archimedean size
$\ll_{a,b,H}(X^HY^{b-1})^k$. Moreover it will be divisible by
$p^{k(k-1)/2}$.  The argument of \cite[\S 3]{Ann} then produces an
auxiliary form $H(\x;\y)$ providing that
\[p\gg_{a,b,H} X^{2H/(k-1)}Y^{2(b-1)/(k-1)}.\]

We now recall that $Y\ll X^a||G||$. Thus on choosing $H$ sufficiently
large we see that it suffices to have $p\gg_{a,b,\ep}(X||G||)^{\ep}X^{2/b}$.
In addition to the form $x_2y_2\partial G/\partial y_1$
that we have already mentioned we now obtain one further form
$H(\x;\y)$ for each $\mathbb{F}_p$-point on the curve $G(s,1;t,1)=0$.
We thus conclude that every solution to $G(\x;\y)=0$ with
$\max|x_i|\le X$ satisfies one of $O_{a,b,\ep}((X||G||)^{\ep}X^{2/b})$
auxiliary conditions $H(\x;\y)=0$.  Here $H$ is a bilinear form
coprime to $G$, with degrees bounded in terms of $a,b$ and $\ep$.  For
each such form, there are $O_{a,b,\ep}(1)$ common solutions to
$G(s,1;t,1)=H(s,1;t,1)=0$ by B\'ezout's Theorem, and $s,t$ determine
$\x,\y$ since these vectors are primitive. This suffices for the lemma.

\section{Further Improvements}

It is possible to reduce slightly the exponent $7/12$ occurring in the
theorem.  Since the improvement is very small we content ourselves
with a very brief sketch of the argument.

The first step is to repeat the analysis of \S 3 taking $K=L=1$ and
obtaining a bi-linear form $F$ in (\ref{bh}), providing that
$M\in[x^{1/2},x]$ satisfies $M\ge (E^{1/3}F^{1/2})^{1+\delta}$.  Note
here that in fact
\[x^{1/2}\le (E^{1/3}F^{1/2})^{1+\delta}\le x\]
for large enough $x$ and small enough $\delta$, since $x^2\ll E^2F\ll
x^2$ and $E,F\gg x^{1/2}$.

When $F$ is bi-linear the Thue equation referred to in \S 4 will have degree
3, and will produce $O_{\eta}(x^{\eta})$ solutions for each interval
$I$, except when $G_3$ is proportional to a cube.  In this case the
corresponding solutions of 
\begin{equation}\label{V}
2x_1x_2y_1+(x_1^2-x_2^2)y_2=1
\end{equation}
lie on a line
\[(x_1,x_2,y_1,y_2)=(x_1^{(0)},x_2^{(0)},y_1^{(0)},y_2^{(0)})+
\lambda(\mu_1,\mu_2,\nu_1,\nu_2)\]
contained in the variety (\ref{V}).

We now write $\cl{M}_0(E,F)$ for the number of quadruples
$(x_1,x_2,y_1,y_2)$ satisfying (\ref{V}), for which 
\[E/2<x_1^2+x_2^2\le E,\;\;\; F/2<y_1^2+y_2^2\le F\]
but which do not lie on a line in the variety (\ref{V}).  We may then
deduce that 
\[\cl{M}_0(E,F)\ll_{\eta}x^{2\eta}E^{1/3}F^{1/2}\ll_{\eta}
x^{2\eta+1-2\psi/3},\] %
where $E=x^{\psi}$ as before. %
Alternatively we can use our previous argument which shows that
\[\cl{M}_0(E,F)\le\cl{M}(E,F)\ll_{\eta}
x^{3\eta+\min(\psi,2-2\psi)/2+9\psi(1-\psi)/8}.\] %
These suffice to show that
\[\cl{M}_0(E,F)\ll_{\eta}x^{3\eta+\varpi}\] %
with %
\[\varpi=\frac{26+\sqrt{433}}{81}=0.57788\ldots\;\; %
(<\frac{7}{12}=0.58333\ldots)\] %
the critical value of $\psi$ being $(55-\sqrt{433})/54$. %

It then remains to consider the form taken by lines lying in the
surface (\ref{V}).  The lines which contain
more than one integral point may be described explicitly, and one is
then able to show that they contribute $O_{\eta}(x^{\eta+1/2})$ to
$\cl{M}(E,F)$.

In this way one may improve the exponent in the theorem to $\varpi$.

The author is grateful to Thomas Reuss for pointing out an error in %
the original version of this final section. %

\bigskip
\bigskip

Mathematical Institute,

24--29, St. Giles',

Oxford

OX1 3LB

UK
\bigskip

{\tt rhb@maths.ox.ac.uk}

\end{document}